\documentclass[12pt]{article}
\usepackage{amsmath}
\usepackage{amsfonts}
\usepackage{amssymb}

\usepackage{amsthm,amsmath,amssymb}
\usepackage{verbatim}
\usepackage{array}
\usepackage[colorlinks=true,citecolor=black,linkcolor=black,urlcolor=blue]{hyperref}

\newtheorem{theorem}{Theorem}

\newtheorem{cor}[theorem]{Corollary}

\newtheorem{definition}[theorem]{Definition}

\theoremstyle{remark}

\usepackage{lipsum}

\begin{document}
\title{{\Large On the size multipartite  Ramsey numbers involving complete graphs}}
\author{\small   Leila Maherani$^{\textrm{b},1}$ \ and \  Maryam Shahsiah$^{\textrm{a},\textrm{b},2}$\\
		\small  $^{\textrm{a}}$ Department of Mathematics, Khansar Campus, University of Isfahan, Isfahan, Iran\\	
	\small  $^{\textrm{b}}$School of Mathematics, Institute for
	Research
	in Fundamental Sciences (IPM),\\
	\small  P.O.Box: 19395-5746, Tehran,
	Iran\\	
	\small \texttt{E-mails: 
		maherani@ipm.ir, m.shahsiah@khc.ui.ac.ir}}
\date {}
\maketitle \footnotetext[1] {\tt This research was in part supported by a grant from IPM (No.1402050046)} \vspace*{-0.5cm}
\footnotetext[2] {\tt This research was in part supported by a grant from IPM (No. 1402050314)} \vspace*{-0.5cm}

%
\setlength{\parindent}{0pt}
\maketitle
\begin{abstract}

 Given two graphs $H$ and $G$, the  size multipartite Ramsey number  $ m_j (H,G )$  is the smallest natural number $t$  such that an arbitrary coloring of the edges
of $K_{j \times t}$, complete multipartite graph whose vertex set is partitioned into $j$ parts each of size $t$,   using  two colors red and blue, necessarily forces a red copy of $H$ or a blue copy of $G$ as a subgraph. The notion of size multipartite Ramsey number  has been introduced by Burger and  Vuuren in 2004. It is worth noting that, this concept  is derived by using the  idea of the original classical Ramsey number, multipartite
Ramsey number and  the size Ramsey number. \\
In this paper, we focus on $ m_j (H,G )$ and  find a  lower bound for $m_j({H,G})$ based on the chromatic number of $H$ and the order of $G$. Also, for graphs $G$ with large maximum degree, we obtain a tight lower bound for $m_j({K_m,G})$. Furthermore, we   determine the order of magnitude of  $	m_j({K_m, K_{1,n}})$, for $j \geq m \geq 3$ and $n\geq 2$. Then we specify the exact values of
$m_j({K_m,K_{1,n}})$ for the cases $m=j$, $m=3$ and $j\equiv 0 \ {\rm or}\ m-2 \ ({\rm mod\ }m-1)$. 
\end{abstract}
\textbf{Keywoeds}: Complete multipartite graphs,  Ramsey number, Multipartite Ramsey number. 

\section{Introduction}
 All graphs $G = (V, E)$ considered in this paper are simple.
 The vertex set and the edge set of graph $G$ are denoted by $V(G)$ and $E(G)$, respectively. We denote by $K_m$  the complete graph on $m$ vertices. The notation  $K_{p_1,p_2, \ldots,p_j}$ represents the
 complete multipartite graph with partite sets of size $p_1,p_2, \ldots,p_j$. When $p_1=p_2= \cdots=p_j=p$, we denote  $K_{p_1,p_2, \ldots,p_j}$ by $K_{j \times p }. $\\

The  \textit{Ramsey number}  $R(H,G)$   is
the smallest integer $m$ such that in every edge coloring of complete graph  ${K}_{m}$  with colors red and blue, there is 
a monochromatic copy of $H$ of color red, or  a monochromatic copy of $G$ of color blue.\\
 
 The multipartite Ramsey number for graphs  has received considerable attention. 
More precisely, it has been interest to see what happens to the Ramsey numbers when
we allow fixed edge deletions from  $K_m$.
 To see more results on multipartite Ramsey numbers we refer the reader to see\cite{Benevides, DGHS,GRSS,GSS}.\\
%
The investigation of the size Ramsey numbers of graphs was initiated by Erd\H{o}s et al. \cite{EFRS} in 1978, by introducing the idea of measuring minimality with respect to size rather than order. Some of the results involving  size Ramsey numbers are \cite{Beck1,Beck2,JM,MOS,RS}.\\

 More recently, using the idea of the original classical Ramsey number, multipartite
 Ramsey number and  size Ramsey number, the notion of size multipartite Ramsey
 number was introduced by Burger and Vuuren, as follow.

 

\begin{definition}{\em \cite{BV2}}\label{Def1}
	Let $j\geq 2$  be an integer number. Given graphs $H$ and $G$, the \textit{size
		multipartite Ramsey number} $ m_j (H,G )$  is the smallest natural number $t$  such that an arbitrary coloring of the edges
	of $K_{j \times t}$ , using  two colors red and blue, necessarily forces a red copy of $H$ or a blue copy of $G$ as a subgraph.
\end{definition}
 Definition \ref{Def1} is a generalization  of the
classical Ramsey numbers in the sense that if $R(K_n,K_s)=\ell$, then $  m_{\ell} (K_{n \times 1}, K_{s \times 1} )=1$.

 In \cite{BV2}, authors  proved that $m_j (K_{n }, K_{s } ) \rightarrow 1$ as $j \rightarrow \infty$.
They also determined exact value of $m_j (K_{n \times l}, K_{s \times p} )$, for some small values of $n,s,l,p$.
\\
 More recently, the size multipartite Ramsey number $m_j (H, G)$ has been studied for special classes
of graphs like paths, cycles, strips and stars. For more details, we refer the reader to see \cite{BV2,  LBS,LSBJ, SS,S,SBU1, SBU2}. It is worth mentioning that, most of these results are for  small cases.\\


In this paper,  we investigate the size multipartite Ramsey number  	$m_j({H,G})$ and  find a  lower bound for $m_j({H,G})$ based on the chromatic number of $H$ and the order of $G$. Also, for graphs $G$ with large maximum degree, we obtain a tight lower bound for $m_j({K_m,G})$. Furthermore, we   determine the order of magnitude of  $	m_j({K_m, K_{1,n}})$ for $j \geq m \geq 3$ and $n\geq 2$, where $K_{1,n}$ is the star with $n+1$ vertices.  We also find the exact values of
$m_j({K_m,K_{1,n}})$ for the cases $m=j$, $m=3$ and $j\equiv 0 \ {\rm or}\ m-2 \ ({\rm mod\ }m-1)$. \\


\subsection{Definitions and preliminaries}

Here, we provide some basic definitions and useful results needed in the following section.\\

Let $G=(V,E)$ be a graph. The {\it degree} $d_G(v)$ of a vertex $v$ is the number  of edges of $G$ incident to $v$.  We denote by $\Delta(G)$ and $\delta(G)$ the  maximum  and minimum degrees of $G$, respectively.  
The {\it chromatic number} of $G$ is the minimum number of colors needed to color the vertices in such a way that no two adjacent vertices have the same color. It is denoted by $\chi(G)$.\\



For given integers $t$ and $j\geq m \geq 2$, let $f(t,j,m)$ be the largest minimum degree $\delta(F) $ among all $F\subseteq K_{j\times t}$ which  contain no copy of $K_{m}$. This notation was introduced by Bollob\'{a}s et al. in  \cite{BES2}.

 \begin{theorem}\label{f} Assume that $j\geq m \geq 3$, $t\geq 1$ are integers. The following results have been proved on $f(t,j,m)$.
 	\begin{itemize}
 			\item [\rm{[i.]}]  {\rm(\cite{BES1})} Let $m=3$.  Then $ f(t,j, 3)= \big\lfloor\frac{j}{2}\big\rfloor t.$ 
 			
 		\item [\rm{[ii.]}]  {\rm(\cite{HS,ST})} If $j$ is even, then we have 
 		$	f(t,j,j)=(j-1)t-\left \lceil \frac{jt}{2(j-1)} \right\rceil$. If $j$ is odd, then  $	f(t,j,j)=f(t,j-1,j-1)+t=(j-1)t-\left \lceil \frac{(j-1)t}{2(j-2)} \right\rceil.$
 		
 				\item [\rm{[iii.]}] {\rm{(\cite{LTZ})} }
 					$\big(j-\big\lceil \frac{j}{m-1}\big\rceil\big)t \leq f(t,j, m)\leq (j-\frac{j}{m-1})t.$
 					
 					\item [\rm{[iv.]}] {\rm{(\cite{LTZ})}}  Let $k\geq 2$,  $m\geq 4$ and  $j=k (m-1)-1.$ Then $ f(t,j, m)=\big(j-\big\lceil \frac{j}{m-1}\big\rceil\big)t.$	
 	\end{itemize}
 \end{theorem}

\vspace*{0.3 cm}
Throughout the paper, 
in an edge coloring of a graph ${G}$ with two colors red and  blue, we use the notations $G_r$ and $G_b$ to refer to  the   subgraphs of ${G}$ that are induced by the edges with colors red and blue, respectively.\\

\indent The  paper is  organized as follows: 
The next section consists of our main results. More precisely, first  we find a lower bound for 	$m_j({H,G})$ based on  $\chi(H)$ and  $n(G)$, Theorem \ref{thm}. Then,  we present a tight lower bound for 	$m_j({K_m,G})$ based on  $\Delta(G)$, Theorem \ref{thm2}. In Theorem \ref{thm1}, we obtain two-sided bounds for $	m_j({K_m, K_{1,n}})$. Next,  in Corollaries \ref{cor1} and \ref{cor2}, we show  these bounds are tight. In Theorem \ref{m=j}, we establish the exact value of 	$	m_j({K_j, K_{1,n}})$ for $j\geq 2$ and finally, we  will close the paper  by determining  the exact value of 	$	m_j({K_3, K_{1,n}})$ in  Theorem \ref{m=3}.
\section {Main results} 

In this section we will investigate the size multipartite Ramsey number for some families of graphs. First we give a lower bound for $m_j({H,G})$ based on  $\chi(H)$ and  $n(G)$.
\begin{theorem}\label{thm}
	Let $j\geq \chi \geq 2$ and $n\geq 2$. If $H$ is a graph with chromatic number $\chi$ and $G$ is a connected graph of order $ n$, then
	\begin{equation*}
		m_j({H,G})\geq  \Big\lfloor \frac{n-1}{\lceil \frac{j}{\chi-1}\rceil} \Big\rfloor+1.
	\end{equation*}
\end{theorem}
{\bf Proof.} Suppose that $t= \big\lfloor \frac{n-1}{\lceil \frac{j}{\chi-1}\rceil} \big\rfloor$ and $F=K_{j\times t}$.  Assume that $V_1, V_2, \ldots, V_j$ are the vertex sets of $F$.  Now, let  $\ell = \big\lceil \frac{j}{\chi-1}\big\rceil$ and  partition  $V_1, V_2, \ldots, V_j$ into $\chi-1$ parts $A_1, A_2, \ldots, A_{\chi-1}$, such that $ |A_i|=\ell$, for each $1\leq i \leq \chi -2$ and $1\leq |A_{\chi -1}|\leq \ell$. 
 Color  all  edges whom end points  belong to the same $A_i$,  $1\leq i \leq \chi-1$, by blue and   the rest by red. It is easy to see that $H\nsubseteq F_r$, since chromatic number of  $F_r$ is at most $\chi -1$. Furthermore, for every component $C$ of $F_b$ we have 
 \begin{align*}
n(C)\leq \ell t= \Big\lceil \frac{j}{\chi-1}\Big\rceil \times \Big\lfloor \frac{n-1}{\lceil \frac{j}{\chi-1}\rceil} \Big\rfloor
 \leq \Big\lceil \frac{j}{\chi-1}\Big\rceil \times \frac{n-1}{\lceil \frac{j}{\chi-1}\rceil} =n-1,
 \end{align*}
where $n(C)$ denotes the number of vertices of $C$.
Therefore, no component of $F_b$ contains $G$. Consequently, $G\nsubseteq F_b$ and $	m_j({H,G})\geq  t+1$.
$\hfill\blacksquare$\\\\

In the following we obtain a tight lower bound for $m_j({K_m, G})$ based on the parameters $m$ and $\Delta(G)$.
\begin{theorem}\label{thm2}
	Let $j\geq m \geq 2$ and $n\geq 2$. Then
	\begin{equation*}
		m_j({K_m, G})\geq  \Big\lfloor \frac{\Delta-1}{\lceil \frac{j}{m-1}\rceil-1} \Big\rfloor+1,
	\end{equation*}
where $\Delta= \Delta(G)$.
\end{theorem}

{\bf Proof.}  First, let $m=2.$  Set $t=\big\lfloor \frac{\Delta-1}{j -1} \big\rfloor$ and $F=K_{j\times t}$. Clearly, $d_v(F)\leq \Delta-1$,  for  every $v\in V(F)$. Consequently,   $m_j({K_2, G})> t.$

Now, let  $j\geq m \geq 3.$  
Suppose that $t'= \big\lfloor \frac{\Delta-1}{\lceil \frac{j}{m-1}\rceil -1} \big\rfloor$ and $F'=K_{j\times t'}$. We consider the edge coloring of $F'$    with two colors red and blue such that $F'_r$ receives  the largest  minimum degree and $K_m\nsubseteq F'_r$. By  Theorem \ref{f}, part [iii],
 	\begin{equation*}
\delta(F'_r)=f(t',j,m) \geq \big(j-\Big\lceil \frac{j}{m-1}\Big\rceil\big)t'. 
	\end{equation*}
	Therefore,  for every vertex $v\in V(F')$, we have 
\begin{align*}
d_{F'_b}(v) &\leq (j-1)t' - \big(j-\Big\lceil \frac{j}{m-1}\Big\rceil\big)t' = \big(\Big\lceil \frac{j}{m-1}\Big\rceil -1\big)t'\\
	&=\big(\Big\lceil \frac{j}{m-1}\Big\rceil -1\big) \times \Big\lfloor \frac{\Delta-1}{\lceil \frac{j}{m-1}\rceil-1} \Big\rfloor \leq \Delta-1.
\end{align*}
 Then, 	$\Delta(F'_b) \leq \Delta-1$ and hence $G \nsubseteq F'_b$. This confirms the desired lower bound.
 $\hfill\blacksquare$\\
\begin{theorem}\label{thm1}
	Let $j\geq m \geq 2$ and $n\geq 2$.  Then
	\begin{align*}
		\Big\lfloor \frac{n-1}{\lceil \frac{j}{m-1}\rceil-1} \Big\rfloor+1\leq	m_j({K_m, K_{1,n}})\leq  \Big\lfloor \frac{n-1}{\frac{j}{m-1}-1} \Big\rfloor+1.
\end{align*}
\end{theorem}
{\bf Proof.} The desired lower bound is concluded by Theorem \ref{thm2}. Therefore, it suffices to prove the claimed upper bound. First, let $m=2$  and $t=  \big\lfloor \frac{n-1}{j-1} \big\rfloor+1$. Suppose that the edges of $F=K_{j\times t}$ are colored with two colors red and blue. If $F_r$ has a $K_2$ as a subgraph, then we are done. Otherwise, for every $v\in V(F)$
 \begin{equation*}
 d_{F_b}(v)= (j-1)\times \big(\Big\lfloor \frac{n-1}{j-1} \Big\rfloor+1\big)> n-1.
  	\end{equation*}
  Therefore, $K_{1,n}\subseteq F_b$ and so $m_j({K_2, K_{1,n}})\leq  \big\lfloor \frac{n-1}{j-1} \big\rfloor+1.$

Now, let $j\geq m \geq 3.$
Assume that $t'=  \big\lfloor \frac{n-1}{\frac{j}{m-1}-1} \big\rfloor+1$ and   the edges of $F'=K_{j\times t'}$ are colored with two colors red and blue. If $F'_r$ has a $K_m$ as a subgraph, then we are done. Otherwise, using Theorem \ref{f}, part [iii],  we have
 \begin{equation*}
  \delta(F'_r)\leq f(t',j,m) \leq \big(j-\frac{j}{m-1}\big)t'.
  	\end{equation*}
   Hence, 
\begin{align*}
	\Delta(F'_b) &\geq (j-1)t' - \big(j-\frac{j}{m-1}\big)t' = \big(\frac{j}{m-1}-1\big)t'\\
	 &=\big(\frac{j}{m-1}-1\big)\times  \Big(\Big\lfloor \frac{n-1}{\frac{j}{m-1}-1} \Big\rfloor+1\Big)\\
	&> \big(\frac{j}{m-1}-1\big)\times   \frac{n-1}{\frac{j}{m-1}-1}=n-1.
\end{align*}
Therefore, $K_{1,n}\subseteq F'_b$ and the proof is completed.
$\hfill\blacksquare$\\  
 
 In the following corollaries, we will show that both of bounds in Theorem \ref{thm1} are tight. 
 
 \begin{cor}\label{cor1}
  Let $k\geq 2$,  $m\geq 4$ and  $j=k (m-1)-1.$  Then
 	$	m_j({K_m, K_{1,n}})=  \big\lfloor \frac{n-1}{\lceil \frac{j}{m-1}\rceil-1} \big\rfloor+1.$
 \end{cor}
	{\bf Proof.} By Theorem \ref{thm1}, it suffices to prove $	m_j({K_m, K_{1,n}})\leq  \big\lfloor \frac{n-1}{\lceil \frac{j}{m-1}\rceil-1} \big\rfloor+1.$
	 Let $t=  \big\lfloor \frac{n-1}{\lceil \frac{j}{m-1}\rceil-1}\big\rfloor+1$. Assume that the edges of $F=K_{j\times t}$ are colored with two colors red and blue. We may suppose that $F_r$ has no copy of $K_m$ as a subgraph, otherwise we are done. Thus,  using Theorem \ref{f}, part [iv],  we have
	  $$\delta(F_r)\leq f(t,j,m) = \big(j-\Big\lceil\frac{j}{m-1}\Big\rceil\big)t.$$ 
	  Hence, 
	\begin{align*}
		\Delta(F_b) &\geq (j-1)t - \big(j-\Big\lceil\frac{j}{m-1}\Big\rceil\big)t = \big(\Big\lceil\frac{j}{m-1}\Big\rceil -1\big)t\\
		&= \big(\Big\lceil\frac{j}{m-1}\Big\rceil -1\big)\times \big(\Big\lfloor \frac{n-1}{\lceil \frac{j}{m-1}\rceil-1}\Big\rfloor+1\big)\\ &>\big(\Big\lceil\frac{j}{m-1}\Big\rceil -1\big)\times  \frac{n-1}{\big\lceil \frac{j}{m-1}\big\rceil-1}=n-1.
	\end{align*}
	Therefore, $K_{1,n}\subseteq F_b$. This completes the proof.
 \begin{cor}\label{cor2}
 	Let $j\geq m \geq  3$  and $m-1$ divides $j.$ Then
 	$	m_j({K_m, K_{1,n}})=  \big\lfloor \frac{n-1}{\frac{j}{m-1}-1} \big\rfloor+1.$
 \end{cor}
 	{\bf Proof.}  Since $\frac{j}{m-1}=\lceil\frac{j}{m-1}\rceil$, the proof is trivial.	
   	$\hfill\blacksquare$\\
\begin{theorem}\label{m=j}
	Let $n \geq 2$ and $j\geq 2$  be positive  integers. If $j$ is even, then 
	$	m_j({K_j, K_{1,n}})=  \big \lfloor \frac{2(n-1)(j-1)}{j} \big\rfloor+1.$
Otherwise, 
$	m_j({K_j, K_{1,n}})=  \big\lfloor \frac{2(n-1)(j-2)}{j-1} \big\rfloor+1.$
\end{theorem}
	
	{\bf Proof. } First let $j$ be even. 
 If $j=2$, then obviously
 $	m_2({K_2, K_{1,n}})=n$ and so we are done. Otherwise, $j\geq 4$. 
Suppose that  $t=  \big\lfloor \frac{2(n-1)(j-1)}{j} \big\rfloor$ and $F= K_{j\times t}$. 
	 Consider a two edge coloring of $F$ with colors red and blue such that $F_r$ receives the largest minimum degree and  $K_j \nsubseteq F_r$.  Using  Theorem \ref{f}, part [ii], we have 
		\begin{equation*}
	\delta(F_r)=f(t,j,j)=  (j-1)t - \Big\lceil\frac{jt}{2(j-1)} \Big\rceil.
	\end{equation*}
Therefore, for every vertex $v \in V(F)$, 
\begin{align*}
	d_{F_b}(v) &\leq \Big\lceil\frac{jt}{2(j-1)} \Big\rceil  = \Big\lceil \frac{j}{2(j-1)}\times \big\lfloor \frac{2(n-1)(j-1)}{j}\big\rfloor \Big\rceil\\
	&\leq \Big\lceil \frac{j}{2(j-1)}\times  \frac{2(n-1)(j-1)}{j} \Big\rceil
	 = n-1.
\end{align*}
Then, $F_b$ contains no  $K_{1,n}$ as a subgraph. Therefore, 	$	m_j({K_j, K_{1,n}}) >t.$ \\
Now, we show that $	m_j({K_j, K_{1,n}}) \leq t+1.$ Consider an arbitrary edge coloring of  $F'=K_{j\times (t+1)}$. We may assume that $K_j \nsubseteq F'_r$. Otherwise, we are done. Therefore, by Theorem \ref{f}, part [ii], we have 
\begin{align*}
	\delta(F'_r) \leq f(t+1,j,j)=(j-1)(t+1)- \Big\lceil\frac{j(t+1)}{2(j-1)}\Big\rceil.
\end{align*}
Thus, 
\begin{align*}
	\Delta(F'_b) &\geq  (j-1)(t+1)-\delta (F'_r) \geq \Big\lceil\frac{j(t+1)}{2(j-1)} \Big\rceil \\ &= \Big\lceil\frac{j}{2(j-1)}\times \big(\big\lfloor \frac{2(n-1)(j-1)}{j}\big\rfloor +1\big) \Big\rceil \\
	 &>  \Big\lceil\frac{j}{2(j-1)}\times  \frac{2(n-1)(j-1)}{j} \Big\rceil  =n-1.
\end{align*}
Thus, there exists a vertex $v$ such that $d_{F_b}(v) \geq n$. Then $K_{1,n}\subseteq F'_{b}$.\\

\indent Now, let $j$ be odd. Assume that $t'=   \big\lfloor \frac{2(n-1)(j-2)}{j-1} \big\rfloor$ and  the edges of $H=K_{j\times t'}$ are  colored with two colors red and blue such that  $H_r$ has the largest minimum degree and   $K_j\nsubseteq H_r$.  Using Theorem \ref{f}, part [ii], we have 
\begin{equation*}
	\delta(H_r)=f(t',j,j)=  (j-1)t' - \Big\lceil\frac{(j-1)t'}{2(j-2)} \Big\rceil.
\end{equation*}
Therefore, for every vertex $v \in V(H)$, 
\begin{align*}
	d_{H_b}(v) &\leq  \Big\lceil\frac{(j-1)t'}{2(j-2)} \Big\rceil  = \Big\lceil\frac{(j-1)}{2(j-2)} \times \big\lfloor \frac{2(n-1)(j-2)}{j-1}\big\rfloor \Big\rceil\\
	&\leq \Big\lceil \frac{(j-1)}{2(j-2)}\times  \frac{2(n-1)(j-2)}{j-1} \Big\rceil
	= n-1.
\end{align*}
Then, $H_b$ contains no  $K_{1,n}$ as a subgraph. Therefore, 	$	m_j({K_j, K_{1,n}}) >t'.$\\
 Now, it suffices to show that $	m_j({K_j, K_{1,n}}) \leq t'+1.$ Consider an arbitrary edge coloring of  $H'=K_{j\times (t'+1)}$ with two colors red and blue. With no  loss of generality,  assume that $K_j \nsubseteq H'_r$.  Otherwise, we are done. By Theorem \ref{f}, part [ii], we have 
\begin{align*}
	\delta(H'_r) \leq f(t'+1,j,j)=(j-1)(t'+1)- \Big\lceil\frac{(j-1)(t'+1)}{2(j-2)}\Big\rceil.
\end{align*}
Therefore, 
\begin{align*}
	\Delta(H'_b) &\geq  (j-1)(t'+1)-\delta (H'_r) \geq \Big\lceil\frac{(j-1)(t'+1)}{2(j-2)}\Big\rceil \\  &= \Big\lceil\frac{(j-1)}{2(j-2)}\times \big(\big\lfloor \frac{2(n-1)(j-2)}{j-1}\big\rfloor +1\big) \Big\rceil \\
	&>  \Big\lceil\frac{j-1}{2(j-2)}\times  \frac{2(n-1)(j-2)}{j-1} \Big\rceil  =n-1.
\end{align*}
This implies that there exists a vertex $v'$ such that $d_{H'_b}(v') \geq n$ and finishes the proof.
$\hfill\blacksquare$\\
\begin{theorem}\label{m=3}
	Let $n \geq 2$ and $j\geq 3$  be positive integers. Then
	\begin{equation*}
		m_j({K_3, K_{1,n}})=  \Big\lfloor \frac{n-1}{{\lceil\frac{j}{2} \rceil} -1}\Big\rfloor+1.
	\end{equation*}
\end{theorem}
	 {\bf Proof. } First, let $t= \Big\lfloor \frac{n-1}{{\lceil\frac{j}{2} \rceil} -1}\Big\rfloor$ and $F=K_{j\times t}$. Assume that the edges of $F$ are  colored with two colors red and blue such that  $F_r$ receives the largest minimum degree and   $K_3 \nsubseteq F_r$. Using Theorem \ref{f}, part [i], we have 
	\begin{align*}
		\delta(F_r)= f(t,j,3)=  \big\lfloor \frac{j}{2}  \big\rfloor t.
	\end{align*}
	Therefore, for every vertex $v \in V(F)$, we have  
	\begin{align*}
	d_{F_b}(v) &\leq (j-1)t-\delta(F_r) = (j-1)t-\big\lfloor \frac{j}{2}  \big\rfloor t\\ &= \big(\big\lceil\frac{j}{2} \big\rceil -1\big)t
		 = \big(\big\lceil\frac{j}{2} \big\rceil -1\big) \times \Big\lfloor \frac{n-1}{\big\lceil\frac{j}{2} \big\rceil -1}\Big\rfloor\\ &\leq \big(\big\lceil\frac{j}{2} \big\rceil -1\big) \times  \frac{n-1}{\big\lceil\frac{j}{2} \big\rceil -1}
		 =n-1.	 	  
	\end{align*}
	Then, this coloring contains no blue copy of $K_{1,n}$. This implies that 	$	m_j({K_3, K_{1,n}}) >t.$ 
	
	Now, we show that $	m_j({K_3, K_{1,n}}) \leq t+1.$ Consider an arbitrary two edge coloring of  $F'=K_{j\times (t+1)}$ with colors red and blue. We may assume that $K_3 \nsubseteq F'_r$. Otherwise, we are done. Therefore, by theorem \ref{f}, Part [i],  
	\begin{align*}
		\delta(F'_r) \leq f(t+1,j,3)=  \big\lfloor \frac{j}{2}  \big\rfloor (t+1).
	\end{align*}
	Consequently, 
	\begin{align*}
		\Delta(F'_b) &\geq  (j-1)(t+1)-\big\lfloor \frac{j}{2}  \big\rfloor (t+1)\\ &= \big(\big\lceil\frac{j}{2} \big\rceil -1\big)(t+1)
		  = \big(\big\lceil\frac{j}{2} \big\rceil -1\big) \times \big(\Big\lfloor \frac{n-1}{\lceil\frac{j}{2} \rceil -1}\Big\rfloor +1\big)\\ 
		&> \big(\big\lceil\frac{j}{2} \big\rceil -1\big) \times \frac{n-1}{\big\lceil\frac{j}{2} \big\rceil -1}
		= n-1.
	\end{align*} 
Therefore,  there exists a vertex  $v$ such that $d_{F'_b}(v) \geq n$.
 It concludes that $K_{1,n} \subseteq F'_b$ and the proof is completed.
$\hfill\blacksquare$\\


\end{document}